%% file: agmg.tex
\documentclass[11pt,a4paper,leqno]{article}

\usepackage{verbatim}
\usepackage{fancyhdr}
\usepackage{graphicx}
\usepackage{geometry}
\usepackage{amssymb}
\usepackage{amsmath}
\usepackage{hyperref}
\usepackage{subfig}
\usepackage{caption}
\usepackage{algpseudocode}
\usepackage{algorithmicx}
\usepackage{algorithm}

\def\usFootnote{Department of Computational and Applied Mathematics, Rice University}
\def\kenFootnote{Stone Ridge Technology}
\title{A GPU Accelerated Aggregation Algebraic Multigrid Method}

\author{R. Gandham\footnote{\usFootnote},\
 K. Elser\footnote{\kenFootnote},\
 Y. Zhang\footnotemark[\value{footnote}]}

\date{}
\begin{document}

\maketitle

\input{abstract}

\input{intro}
\input{aggregationAMG}
\input{implementation}
\input{experiments}
\input{conclusions}

\clearpage

\bibliographystyle{siam}
\bibliography{refs}

\pagebreak
\end{document}

%% file: abstract.tex
\begin{abstract}
We present an efficient, robust and fully GPU-accelerated aggregation-based algebraic multigrid preconditioning technique for the solution of large sparse linear systems. These linear systems  arise from the discretization of elliptic PDEs. The method involves two stages, setup and solve. In the setup stage, hierarchical coarse grids are constructed through aggregation of the fine grid nodes. These aggregations are obtained using a set of maximal independent nodes from the fine grid nodes. We use a ``fine-grain'' parallel algorithm for finding a maximal independent set  from a graph of strong negative connections. The aggregations are combined with a piece-wise constant (unsmooth) interpolation from the coarse grid solution to the fine grid solution, ensuring  low setup and interpolation cost. The grid independent convergence is achieved by using recursive Krylov iterations (K-cycles) in the solve stage. An efficient combination of  K-cycles and standard multigrid V-cycles is used as  the preconditioner for Krylov iterative solvers such as generalized minimal residual and conjugate gradient. We compare the solver performance with other solvers based on smooth aggregation and classical algebraic multigrid methods.
\end{abstract}


%% file: intro.tex
\section{Introduction}
In the numerical simulation of various physical phenomena, discretizations of partial differential equations result in very large sparse linear system of equations. Computing solutions of these linear systems is a crucial component of the overall simulation. Consequently  there has been increasing interest in developing fast linear solvers over the last two decades. Traditionally, these linear solvers are accelerated by using parallel programming techniques on clusters of CPUs. However, hardware architectures such as the one found in the  graphics processing unit (GPU) provide efficient, inexpensive alternatives for these computations. The recent developments in general purpose GPUs (GPGPUs) for scientific computations make it feasible to accelerate the numerical simulations by an order of magnitude. This can be achieved by significantly modifying the traditional algorithms to suit the massively parallel hardware architecture of GPUs. \\

\noindent Multigrid methods are among the most efficient and popular solution techniques for solving the linear systems arising from the discretization of elliptic PDEs. Multigrid methods typically fall in to one of two classes: geometric multigrid and algebraic multigrid (AMG). Geometric multigrid methods require  prior knowledge of the underlying discretization and grid hierarchy, where as algebraic multigrid methods only require the entries of matrix. Furthermore, the algebraic multigrid  coarsening process automatically ensures a reduction of the error components that are not reduced by smoothing \cite{ruge1987algebraic}. Algebraic multigrid requires little or no tuning for various applications and hence make it very easy to use in the form of an efficient ``black-box''  solver. \\

\noindent Algebraic multigrid methods involve the construction of  a hierarchy of matrices (or grids) using the entries of the matrix. At each level in the hierarchy, smoothing is performed to remove the high frequency errors.  The low frequency errors are removed by recursively restricting the residual to coarse grids and performing smoothing on the coarse grids. The coarse grid solutions are then interpolated from the coarse grid to a fine grid using an interpolation operator. \\

\noindent Classical algebraic multigrid methods obtain the hierarchical grids   by partitioning the nodes into coarse and fine grid nodes. The coarse grid nodes form a coarse level, and an interpolation operator is defined via a weighted sum of the coarse grid nodes. The restriction operator in general is the transpose of the interpolation operator \cite{stuben2001review}. \\

\noindent In contrast, aggregation algebraic multigrid methods obtain the hierarchical grids by aggregating few fine grid nodes to form a coarse grid node (see Figure \ref{fig:example_graph}). An interpolation is defined via  piecewise constant interpolation from coarse grid node  to a fine grid node. This results in very sparse interpolation and restriction operations compared to classical AMG. The interpolation matrix has exactly one nonzero entry per row, which reduces memory requirements and improves efficiency of the interpolation operation. However, these aggregation schemes are not popular since they do not provide grid independent convergence. The convergence is often improved by smooth interpolation (``\emph{smooth aggregation}'') \cite{vanek1996algebraic}, but  smooth aggregation typically produces a hierarchy with  very dense coarse level matrices (see Figure \ref{fig:compare_sparsity}), leading to expensive matrix-vector product computations and increasing the memory requirements. Furthermore, unlike classical algebraic multigrid, smooth aggregation multigrid is not robust for various applications. \\

\noindent The opposite view point to smooth aggregation, presented in \cite{notay2010aggregation} uses unsmooth aggregation along with K-cycles as preconditioning for iterative Krylov methods, provide an efficient alternative to smooth aggregation schemes. In this work, we accelerate these unsmooth aggregation methods using GPUs and discuss the algorithmic changes considered for  the acceleration, both mathematical and computational. We also discuss the advantages in using these schemes  where  we can efficiently reuse the sparsity structure of the hierarchical matrices. \\

\noindent Some of the initial works for aggregation AMG focused on accelerating only the solve phase \cite{emans2012steps}. As a consequence, the setup phase becomes a bottleneck due to the serial nature of algorithms \cite{notay2010aggregation} for forming the aggregates. There has been development in accelerating the setup phase for classical AMG \cite{esler2012gampack} and for smooth aggregation \cite{bell2012exposing}. We build on these ideas and develop fully accelerated unsmooth aggregation AMG as a part of \verb=GAMPACK=, a GPU accelerated algebraic multigrid package. \\

\noindent For forming the aggregates on GPUs, we follow a fine-grain parallel maximal independent set (MIS) algorithm proposed in \cite{bell2012exposing} with some modifications described in algorithm (\ref{alg:pmis2}). With this algorithm on a GPU, we see that the setup is almost three times faster than classical AMG\ setup.\\

\noindent In Section \ref{sec:methods}, we describe the algorithms for  the setup and solution phase of unsmooth aggregation AMG. In Section \ref{sec:implementation}, we briefly explain GPU acceleration of the algorithms described. In Section \ref{sec:results}, we provide numerical results for various classes of problems and compare the performance with  \verb=GAMPACK= classical AMG and \verb=CUSP= smooth aggregation. 

%% file: aggregationAMG.tex
\section{Aggregation AMG}
\label{sec:methods}
Aggregation AMG is used as a preconditioner for iterative Krylov solvers such as  generalized minimal residual (GMRES) and conjugate gradient (CG). Consider a linear system of equations,
\begin{equation}
\label{eq:axb}
A x = b,
\end{equation}
where $A \in \mathbb{R}^{N \times N}$ is a coefficient matrix, $b \in \mathbb{R}^{N \times 1}$ is a right hand side vector, $x \in \mathbb{R}^{N \times 1}$  is the solution vector and $N$  is the number of unknowns. 
The coarse grid matrices are defined recursively, consider that $A_0 = A$ and $n_{k+1} < n_k, \, n_0 = N$. The coarse grid operators ($A_{k+1}$) are defined by,
\begin{equation}
A_{k+1_{}} = R_{k}  A _{k} P_{k},  \qquad k \, = 0, 1, .. ,L
\end{equation}
where $P_k \in \mathbb{R}^{n_k \times n_{k+1}}$ is an interpolation operator, $R_k$ ($=P_k ^T$) is a restriction operator, $A_k \in \mathbb{R}^{n_k \times n_k}$ is a fine grid coefficient matrix, while $A_{k+1} \in \mathbb{R} ^{n_{k+1} \times n_{k+1}}$ is a coarse grid coefficient matrix. The coarse grid operators are obtained until the size of the coarse grid operator ($n_L$) is sufficiently small, to obtain the exact solution using a direct solver in a reasonable time. \\

\noindent Aggregation algebraic multigrid uses the prior knowledge of  null space or near null space of the linear operator that corresponds to low-energy error. In our work, we consider only one near null space vector of the linear operator and propagate the null space to  hierarchical operators during the setup stage. For level $k$, the near null space vector of matrix $A_k$ is denoted by $B_k$.

\subsection{Setup phase}
Due to simple structure of interpolation matrix, unsmooth aggregation schemes have the advantage of low setup and interpolation costs compared to that of smooth aggregation schemes or classical algebraic multigrid schemes.  The algorithm for the setup is outlined in Algorithm (\ref{alg:setup}).
\begin{algorithm}[!h]
\caption{Aggregation AMG setup}
\label{alg:setup}
\begin{algorithmic}[1]
\State \textbf{Input\,: } $A_0$, $B_0$
\State \textbf{Output: } $A_k$, $B_k$, $P_k$, $R_k$ for $k=1,2,..L$
\For {$k = 0, 1, 2,..$ until $n_L $ is small}
\State $C_k \leftarrow $ \verb=strength=($A_k$)
\State $Agg_k \leftarrow $ \verb=aggregate=($C_k$)
\State $P_k, B_{k+1} \leftarrow $ \verb=prolongate=($Agg_k$, $B_k$)
\State $R_k \leftarrow P_k ^{T}$
\State $A_{k+1} \leftarrow R_k A_k P_k$
\EndFor
\end{algorithmic}
\end{algorithm}

\subsubsection{Strong connections}
The aggregation scheme groups a set of fine nodes that are strong negatively coupled with each other into a coarse node, called an aggregate. In general, the strength of the connection between two nodes is defined based on the matrix coefficient corresponding to the connection. Two of the popular strength of measures that are shown to be robust are symmetric strength \cite{vanek1996algebraic} and classic strength \cite{ruge1987algebraic}. We observe that classic strength ensures a consistent ratio of the number of aggregates to the fine grid nodes for a large class of problems. A graph of strong connections ($C$) is constructed based on classic strength of connection given by,
\begin{equation}
C_{ij} =
\begin{cases} 1 & \text{if  } \, -s_i A_{ij} > \alpha \,  \displaystyle\max_{-s_i A_{ij} > 0}  -s_i A_{ij},  \qquad s_i = \text{sign}(A_{ii})\\
0 & \text{otherwise}
\end{cases}
\end{equation}
Here $\alpha$ is a threshold parameter for strong and weak connections. The coarsening ratio can be tuned by adjusting $\alpha$. We choose $\alpha = 0.5$ as default for systems arising from 3D problems, and $\alpha = 0.25 $ for systems arising from 2D problems. This selection of $\alpha$ ensures a good ratio of fine grid nodes to coarse grid nodes for a large class of problems that we considered in this paper.
\subsubsection{Construction of aggregates}
The pairwise aggregation scheme proposed in \cite{notay2006aggregation} is  simple and efficient on CPUs, but is not readily parallelizable on GPUs. In addition, this scheme  requires two passes of pairwise aggregation to obtain a sufficient reduction in the number of unknowns for the coarse grid, increasing the setup cost. The aggregation algorithm based on maximal independent set, proposed in \cite{bell2012exposing} is highly parallelizable on GPUs and constructs the aggregates faster compared to pair-wise algorithms. In the aggregation scheme, a set of nodes that are maximally independent are selected as root nodes, and aggregates are formed by grouping each of these nodes with their neighbors. \\

\noindent An \emph{independent set} is a set of nodes, in which no two of them are adjacent and it is \emph{maximal independent set} (MIS) if it is not a subset of any other independent set. The generalization of MIS is MIS($k$), in which the distance between any two independent nodes is greater than $k$ and for every other node there is at least one independent node that is within distance less than or equal to $ k$. For the construction of aggregates, we use MIS($2$) nodes as root nodes. For more aggressive coarsening, MIS($k$), $k>2$ can be used, but these may result in poor interpolation leading to slower convergence. \\

\begin{algorithm}
\caption{Parallel Maximal Independent Set (2) }
\label{alg:pmis2}
\begin{algorithmic}[1]
\algnotext{EndFor}
\algnotext{EndIf}
\algnotext{EndWhile}
\State \textbf{Input:} $C$, $N \times N$ sparse matrix;
\State \textbf{Output :} $s$,  set of MIS(2) nodes

\State $s \leftarrow \{ 0, 0, 0, ..., 0\}$ \hfill \{initialize the state as undecided \}
\State $r \leftarrow $ random \hfill \{ generate random numbers $\in (0, 1)$ \}

\ForAll {$i \in I$}
\State $v_i = \# \{j : C_{ji} = 1 \}  + r_i$ \hfill \{ number of strong influencing connections\}
\EndFor


\While {$\{ i \in I: s_i = 0 \} \neq \emptyset $}

\ForAll {$i \in I$}
\State $T_{i} \leftarrow (s_i, v_i, i)$  \hfill \{ initialize tuples \}
\EndFor


\For {$d = 1, 2$}
\ForAll {$i \in I$} \hfill \{for each node in parallel \}
\State $t \leftarrow T_i$
\For {$j \in N_i$}
\State $t \leftarrow \max(t, T_j)$ \hfill \{compare with the strong neighbours \}
\EndFor
\State $\hat{T}_i \leftarrow t$
\EndFor
\State $T = \hat{T}$
\EndFor


\For {$i \in I$} \hfill \{ for each node in parallel \}
\State $(s_{\max}, v_{\max}, i_{\max}) \leftarrow  T_i$
 \If {$s_i = 0$} \hfill \{ if undecided \}
 \If {$i_{\max} = i$} \hfill \{ if maximal \}
 \State $s_i \leftarrow 1$ \hfill \{ mark as MIS \}
 \Else \If {$s_{\max} = 1$} \hfill \{ else..\}
 \State $s_i \leftarrow -1$ \hfill \{ mark as non MIS \}
\EndIf
\EndIf
\EndIf
\EndFor

\EndWhile
\end{algorithmic}
\end{algorithm}

\noindent In order to aggregate only strongly connected nodes the algorithm utilizes  the graph of matrix $C$  instead of $A$.  The parallel MIS Algorithm (\ref{alg:pmis2}) is similar to the proposed algorithm in \cite{bell2012exposing}, additionally it incorporates the number of strongly influencing connections for robustness \cite{notay2010aggregation}. Node $j$ is a strongly influencing node of $i$ if $A_{ji}$ is a strong connection, in other words $C_{ji} = 1$. After the MIS(2) nodes are obtained, every other node $i$, is aggregated with its nearest MIS node.
\subsubsection{Interpolation}
The interpolation operator for projecting  errors from a coarse grid to fine grid, is defined from the aggregations  and the near null space vectors ($B_k$) of the linear operator ($A_k$). In this paper we consider only one near null space vector for constructing the interpolation operator. This results in exactly one non zero entry per row in the interpolation matrix and is given by,
\begin{equation}
\label{eq:null_space}
B_k = P_k B_{k+1}, \qquad P_k ^T P_k = I.
\end{equation}
Equations (\ref{eq:null_space}) ensure the near null space vector ($B_k$) to be in the range space of ($P_k$) and the interpolation matrix is orthonormal. This is done by copying the entries of $B_k$ to the sparsity pattern of the interpolation matrix and by normalizing the columns.  The restriction operator carries the residual from the fine grid to the next coarse grid in the hierarchy. A  parallel algorithm for computing the interpolation operator is presented in Algorithm (\ref{alg:interp}).

\begin{algorithm}
\caption{Interpolation}
\label{alg:interp}
\begin{algorithmic}[1]
\algnotext{EndFor}
\State \textbf{Input:} $Agg$, $B_k$
\State \textbf{Output:} $P_k$, $R_k$

\ForAll {$1 \le i \le n_k $} \hfill \{ for each fine node in parallel \}
\State $P_k (i,\, I) \leftarrow B_k(i), \qquad I = Agg(i)$ \hfill \{ copy null space vector entries \}
\EndFor

\State $R_k \leftarrow P_k ^{T}$ \hfill \{ transpose \}

\ForAll {$ 1 \le i \le n_{k+1}$} \hfill \{ for each coarse node in parallel \}
\State $B_{k+1}(i) \leftarrow \lVert R_k(i, \, :) \rVert _{2}$ \hfill \{ $L_2$ norm of each row \}
\State $R_k(i,\, :) \leftarrow R_k(i, \, :)/B_{k+1}(i)$ \hfill \{ normalize row \}
\EndFor

\ForAll {$ 1 \le i \le n_k $} \hfill \{ for each fine node in parallel \}
\State $P_k(i, \, I) \leftarrow P_k(i,\, I)/B_{k+1}(I)$ \hfill \{ normalize interpolator \}
\EndFor

\end{algorithmic}
\end{algorithm}
\subsection{Solve phase}
We use a multigrid cycle as a preconditioner in a Krylov iteration which is based on either CG or GMRES. Both of these methods can be used for symmetric matrices while only GMRES\ can be used for non-symmetric matrices. \\

\noindent For unsmooth aggregation, standard V-cycle multigrid given in Algorithm (\ref{alg:vcycle}) does not provide grid independent convergence. In contrast, the K-cycle Algorithm (\ref{alg:kcycle}), presented in \cite{notay2010aggregation}, addresses this issue by recursively applying Krylov iterations on the coarse levels. \\

\begin{algorithm}[!h]
\caption{ V-cycle AMG}
\label{alg:vcycle}
\begin{algorithmic}[1]
\algnotext{EndFor}
\algnotext{EndIf}
\State $x_k \leftarrow$ \verb=Vcycle= $(k, b_k, x_k)$
\State \textbf{Input  :} level $k$, rhs $b_k$ and initial guess $x_k$
\State \textbf{Output :} updated solution $x_k$
\State $x_k \leftarrow S_k(b_k, A_k, x_k)$ \hfill \{ pre-smoothing \}
\State $r_k \leftarrow b_k - A_k x_k$ \hfill \{ compute residual \}
\State $r_{k+1} \leftarrow R_k r_k$ \hfill \{ restrict the residual to coarse-grid \}
\State
\If  {$k+1 = L$}
\State $x_{k+1} \leftarrow A_{k+1} ^{-1} r_{k+1}$ \hfill \{ exact solution of coarse-grid \}
\Else
\State $x_{k+1} \leftarrow$ \verb=Vcycle= $(k+1, r_{k+1}, 0)$ \hfill \{ recursion \}
\EndIf
\State
\State $x_k \leftarrow x_k + P_k x_{k+1}$ \hfill \{ prolongation \}
\State $x_k \leftarrow S_k(b_k, A_k, x_k)$ \hfill \{ post-smoothing \}
\end{algorithmic}
\end{algorithm}

\begin{algorithm}[!h]
\caption{K-cycle AMG}
\label{alg:kcycle}
\begin{algorithmic}[1]
\algnotext{EndFor}
\algnotext{EndIf}
\State $x_k \leftarrow$ \verb=Kcycle= $(k, b_k, x_k)$
\State \textbf{Input  : } level $k$, rhs $b_k$ and initial guess $x_k$
\State \textbf{Output : } updated solution $x_k$
\State $x_k \leftarrow S_k(b_k, A_k, x_k)$ \hfill \{ pre-smoothing \}
\State $r_k \leftarrow b_k - A_k x_k$ \hfill \{ compute residual \}
\State $r_{k+1} \leftarrow R_k r_k$ \hfill \{ restrict the residual to coarse-grid \}
\State
\If {$k+1 = L$}
\State $x_{k+1} \leftarrow A_{k+1} ^{-1} r_{k+1}$ \hfill  \{ exact solution of coarse-grid \}
\Else
\State $c_{k+1} \leftarrow$ \verb=Kcycle= $(k+1, r_{k+1}, x_{k+1})$ \hfill \{ inner first Krylov iteration \}
\State $v_{k+1} \leftarrow A_{k+1} c_{k+1}$
\State $\rho_{1} \leftarrow c_{k+1} ^{T} v_{k+1}, \qquad \alpha_{1} \leftarrow c_{k+1}^{T} r_{k+1} \qquad $ if CG
\State $\rho_{1} \leftarrow \lVert v_{k+1} \rVert ^{2}, \qquad \alpha_{1} \leftarrow v_{k+1}^{T} r_{k+1} \qquad $     if GMRES

\State
\State $\tilde{r}_{k+1} \leftarrow r_{k+1} - \frac{\alpha_{1}}{\rho_{1}} v_{k+1}$
\If {$\lVert \tilde{r}_{k+1}  \rVert \le t \lVert r_{k+1} \rVert$}
\State  $x_{k+1} \leftarrow \frac{\alpha_{1}}{\rho_{1}}c_{k+1}$
\Else
\State $d_{k+1} \leftarrow$  \verb=Kcycle= $(k+1, \tilde{r}_{k+1}, x_{k+1} )$ \hfill \{ inner second Krylov iteration \}
\State  $w_{k+1} \leftarrow A_{k+1}d_{k+1}$
\State $\gamma \leftarrow d_{k+1} ^{T} v_{k+1}, \qquad \beta \leftarrow d_{k+1} ^{T} w_{k+1},  \qquad \alpha_2 \leftarrow d_{k+1} ^{T} \tilde{r}_{k+1} \qquad $ if CG
\State $\gamma \leftarrow w_{k+1} ^{T} v_{k+1}, \qquad \beta \leftarrow \lVert w_{k+1} \rVert ^{2}, \qquad \alpha_2 \leftarrow w_{k+1} ^{T} \tilde{r}_{k+1} \qquad $ if GMRES
\State $\rho_2 \leftarrow \beta - \frac{\gamma ^2}{ \rho_1}$
\State $x_{k+1} \leftarrow (\frac{\alpha_1}{\rho_1} - \frac{\gamma \alpha_2}{\rho_1 \rho_2}) c_{k+1} + \frac{\alpha_2}{\rho_2}d_{k+1}$
\EndIf
\EndIf
\State
\State $x_k \leftarrow x_k + P_{k}x_{k+1}$ \hfill \{ prolongation \}
\State $r_k \leftarrow \tilde{r}_k - A_k x_k$ \hfill \{ compute new residual \}
\State $x_k \leftarrow S_k (r_k, A_k, x_k)$ \hfill \{ post-smoothing \}
\end{algorithmic}
\end{algorithm}

\noindent With K-cycles, the number of coarse grid corrections can grow exponentially, leading to a large number of iterations on the coarser grids. Computing these coarse grid corrections on GPUs is inefficient because of the smaller number of nodes. Interestingly, the experimental results suggest that it is sufficient to use K-cycles for only  two levels at the top and V-cycles for the remaining levels to achieve grid independence convergence. This improves the overall runtime performance of   the aggregation AMG on GPUs  even though the number of iterations increases slightly.

\subsubsection{Smoothing}
In the above algorithms describing V- and K-cycles, $S_k$ represents a smoothing operator that removes high frequency errors corresponding to the matrix $A_k$. We use one pre- and one post- smoothing step for each level. We consider Jacobi, damped Jacobi and symmetric Gauss-Seidel smoothing, all of which are symmetric. Jacobi and damped Jacobi are easily parallelizable on CPUs and GPUs. However, the Gauss-Seidel algorithm is completely serial in nature and cannot be ported to GPUs efficiently. A multi-color variant of symmetric Gauss Seidel algorithm has shown to scale well on GPUs for large matrices with fewer numbers of colors. This modified\  algorithm adds a setup cost of the graph coloring and a reordering of the matrix, and also results in slow down of convergence.  Numerical results indicate that damped Jacobi is both scalable and efficient, and requires low setup cost compared to multi-color Gauss Seidel. In addition, damped Jacobi produces solutions that are independent of number of GPUs/processors unlike multi-color Gauss Seidel. \\

\noindent Damped Jacobi iteration is given by,
\begin{equation}
\label{eq:damped_jacobi}
x \leftarrow x + \omega D^{-1} (b - Ax), \qquad \omega = \frac{4}{3} \frac{1}{\rho (D^{-1} A)}\, ,
\end{equation}
where $\rho (D^{-1} A)$ is  the spectral radius of the matrix $D^{-1}A$ and it is estimated from the eigenvalues of the Hessenberg matrix obtained using $m \, (\le 5)$ Arnoldi iterations.

%% file: implementation.tex
\section{Implementation}
\label{sec:implementation}
We use the CUDA programming model for our GPU implementations and the OpenMP threading model for CPU counter parts. For the setup phase, the matrices are stored in compressed sparse row (CSR) format for both CPU and GPU implementations. For sparse matrix matrix products (\verb=spmm=),   we also use the coordinate (COO)\ format for storing one of the multipliers.      \\

\noindent In the solve phase, matrix vector multiplication \verb=axpy=, is the dominant component. For the CPU implementations we store the matrices in CSR format and for GPU implementations we use a hybrid matrix format \cite{bell2008efficient} that is a combination of the ELL and COO formats. For multi-GPU \verb=axpy=, the matrix entries corresponding to the columns that are processed by other processors are stored in COO format. This way, the multiplication with the matrix stored in ELL,  and the communication of the vector data can be overlapped. We use CUDA peer-to-peer transfers for inter GPU communication of vector entries.  \\

\noindent For parallel primitive operations such as reduction, sort by key, scan, and gather required for the aggregation algorithm, we use the libraries \verb=thrust= \cite{Thrust} and/or \verb=GAMPACK=. \\

\noindent For the coarsest level, to solve the linear system exactly, we use the CPU based direct solver \verb=SuperLU= \cite{superlu99}, for both GPU and CPU implementations. Only one CPU thread is used for solving the coarse linear system. We ensure that dimension of the coarse linear system is about $600$ for GPU implementations to minimize the overall runtime.

\subsection{Galerkin product}
After the interpolation and restriction matrices are computed at a level $k$, the coarse grid operator ($A_{k+1} = R_k A_k P_k$) is constructed using two matrix-matrix multiplications,  $A_k \times P_k$ and $R_k \times (A_k \times P_k)$. All of these matrices are stored in CSR format. A GPU implementation of  sparse matrix matrix multiplication is described in \cite{esler2012gampack}. We  describe another approach that will be useful in the solution of nonlinear systems of equations, where the aggregations are not  constructed at every Newton step when  the change in matrix entries is not significant. In these cases, sparsity pattern of the matrices in the hierarchy do   not change since the pattern depends only on the aggregations and the sparsity pattern  of $A_0$. This is true only if the sparsity pattern of $A_0$ does not change at every Newton step.  \\

\noindent The Galerkin product can be represented as,
\small
\begin{equation}
(A_{k+1})_{IJ} = \sum_{i \in G_I} \sum_{j \in G_J} (R_k) _{Ii} (A_k) _{ij} (P_k)_{jJ} = \sum_{i \in G_I} \sum_{j \in G_J} (P_k) _{iI} (A_k) _{ij} (P_k)_{jJ}
\end{equation}
\normalsize
Where, $G_I$ is a set of nodes that form the $I^{th}$ aggregate. Each of the fine nodes is grouped into exactly one aggregate. Therefore, for each nonzero of the fine grid matrix $A_{ij} \neq 0$, the index of corresponding coarse grid matrix is $IJ$,  where $I$ and $J$ are aggregation indices corresponding  to $i$ and $j$ respectively. These coarse grid nonzero row and column indices are tabulated and sorted based on lexicographic ordering of ordered pair $(I,J)$. The number of nonzeroes in coarse grid matrix is obtained by counting the unique indices $(I,J)$. All the coefficients corresponding to a same index, will sum to a unique nonzero of the coarse grid. This summation is done using segmented reduction. The pattern of this computation can be repeated by storing the sorting indices and the segmented reduction indices. This way, the entire Galerkin product can be performed in less than $1/10^{th}$ of the regular \verb=spmm= compute time. Note that the interpolation operator is not reconstructed, since its  entries  depend on the near null space vector and not directly on the entries of the level matrix. For the cases in which we do not reuse the hierarchy information; we use GPU sparse matrix matrix multiplications.\\

\begin{figure}
\begin{center}
\includegraphics[trim=0cm 18cm 2cm 0cm,clip=true,width=0.5\textwidth]{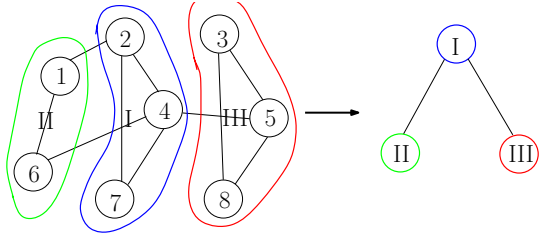}
\caption{\emph{Example of aggregation. Roman numerals indicate the aggregation index.}}
\label{fig:example_graph}
\end{center}
\end{figure}

\noindent Consider a matrix $A_k$, with sparsity graph shown in Figure (\ref{fig:example_graph}). The aggregates are shown in different colors and enumerated in Roman numerals. The sequence of operations to compute the coarse grid matrix is illustrated in the Equations (\ref{eq:galerkin_prod}). Initially the nonzero coefficients $A_{ij}$'s are sorted based on the corresponding aggregation indices $(I,J)$ to obtain unique $(I,J) $ pairs.  The Galerkin product is then finished by adding all the duplicate pairs $(I,J)$ to corresponding nonzero of the coarse grid matrix $A_{k+1}$. Note that each nonzero entry $(A_{k})_{ij}$ is multiplied with $(P_k)_{jJ}$ and $(R_k)_{Ii}$ before sorting the coefficients. \\

\small
\begin{equation}
\label{eq:galerkin_prod}
\begin{bmatrix}i & j & (A_{k})_{ij} \\
1 & 1 & a_{11} \\
1 & 2 & a_{12} \\
1 & 6 & a_{16} \\
2 & 1 & a_{21} \\
2 & 2 & a_{22} \\
2 & 4 & a_{24} \\
2 & 7 & a_{27} \\
3 & 3 & a_{33} \\
3 & 5 & a_{35} \\
3 & 8 & a_{38} \\
4 & 2 & a_{42} \\
4 & 4 & a_{44} \\
4 & 5 & a_{45} \\
4 & 6 & a_{46} \\
4 & 7 & a_{47} \\
5 & 3 & a_{53} \\
5 & 4 & a_{54} \\
5 & 5 & a_{55} \\
5 & 8 & a_{58} \\
6 & 1 & a_{61} \\
6 & 4 & a_{64} \\
6 & 6 & a_{66} \\
7 & 2 & a_{72} \\
7 & 4 & a_{74} \\
7 & 7 & a_{77} \\
8 & 3 & a_{83} \\
8 & 5 & a_{85} \\
8 & 8 & a_{88} \\
\end{bmatrix} \xrightarrow{(i,j)\to (I,J)}
\begin{bmatrix}I & J & (A_{k})_{ij}\\
2 & 2 & a_{11} \\
2 & 1 & a_{12} \\
2 & 2 & a_{16} \\
1 & 2 & a_{21} \\
1 & 1 & a_{22} \\
1 & 1 & a_{24} \\
1 & 1 & a_{27} \\
3 & 3 & a_{33} \\
3 & 3 & a_{35} \\
3 & 3 & a_{38} \\
1 & 1 & a_{42} \\
1 & 1 & a_{44} \\
1 & 3 & a_{45} \\
1 & 2 & a_{46} \\
1 & 1 & a_{47} \\
3 & 3 & a_{53} \\
3 & 1 & a_{54} \\
3 & 3 & a_{55} \\
3 & 3 & a_{58} \\
2 & 2 & a_{61} \\
2 & 1 & a_{64} \\
2 & 2 & a_{66} \\
1 & 1 & a_{72} \\
1 & 1 & a_{74} \\
1 & 1 & a_{77} \\
3 & 3 & a_{83} \\
3 & 3 & a_{85} \\
3 & 3 & a_{88} \\
\end{bmatrix}  \xrightarrow{sort}
\begin{bmatrix}I & J & \\
1 & 1 & a_{22} \\
1 & 1 & a_{24} \\
1 & 1 & a_{27} \\
1 & 1 & a_{42} \\
1 & 1 & a_{44} \\
1 & 1 & a_{47} \\
1 & 1 & a_{72} \\
1 & 1 & a_{74} \\
1 & 1 & a_{77} \\
1 & 2 & a_{21} \\
1 & 2 & a_{46} \\
1 & 3 & a_{45} \\
2 & 1 & a_{12} \\
2 & 1 & a_{64} \\
2 & 2 & a_{11} \\
2 & 2 & a_{16} \\
2 & 2 & a_{61} \\
2 & 2 & a_{66} \\
3 & 1 & a_{54} \\
3 & 3 & a_{33} \\
3 & 3 & a_{35} \\
3 & 3 & a_{38} \\
3 & 3 & a_{53} \\
3 & 3 & a_{55} \\
3 & 3 & a_{58} \\
3 & 3 & a_{83} \\
3 & 3 & a_{85} \\
3 & 3 & a_{88} \\
\end{bmatrix}  \xrightarrow{reduce}
\begin{bmatrix}I & J & (A_{k+1})_{IJ} \\
1 & 1 & a_{22} + a_{24}  \\
& & +a_{27} +\ a_{42} \\
& &+a_{44}+a_{47}\\
&&+a_{72}+a_{74} \\
&&+a_{77} \\
1 & 2 & a_{21}+a_{46} \\
1 & 3 & a_{45}\\
2 & 1 & a_{12}+a_{64} \\
2 & 2 & a_{11}+a_{16} \\
& &+a_{61}+a_{66} \\
3 & 1 & a_{54}\\
3 & 3 & a_{33}+a_{35} \\
&&+a_{38}+a_{53}\\
&&+a_{55}+a_{58}\\
&&+a_{83}+a_{85}\\
&&+a_{88} \\
\end{bmatrix}
\end{equation}
\normalsize

%% file: experiments.tex
\section{Numerical Experiments}
\label{sec:results}

In order to determine the efficiency of the aggregation AMG solver, we compare its performance with that of \verb=GAMPACK= implementation of classical AMG solver. This classical AMG solver was compared with the well known classical AMG solver \verb=HYPRE= \cite{esler2012gampack}. We compare the best performances of both  solvers for each case. We used a flexible GMRES as the outer Krylov solver in all test cases. \\

\begin{figure}
\begin{center}
\subfloat[Smooth, level = 1]{
\begin{minipage}[c]{0.25\linewidth}
\centering
\includegraphics[trim=4cm 6cm 4cm 6cm,clip=true,width=\textwidth]{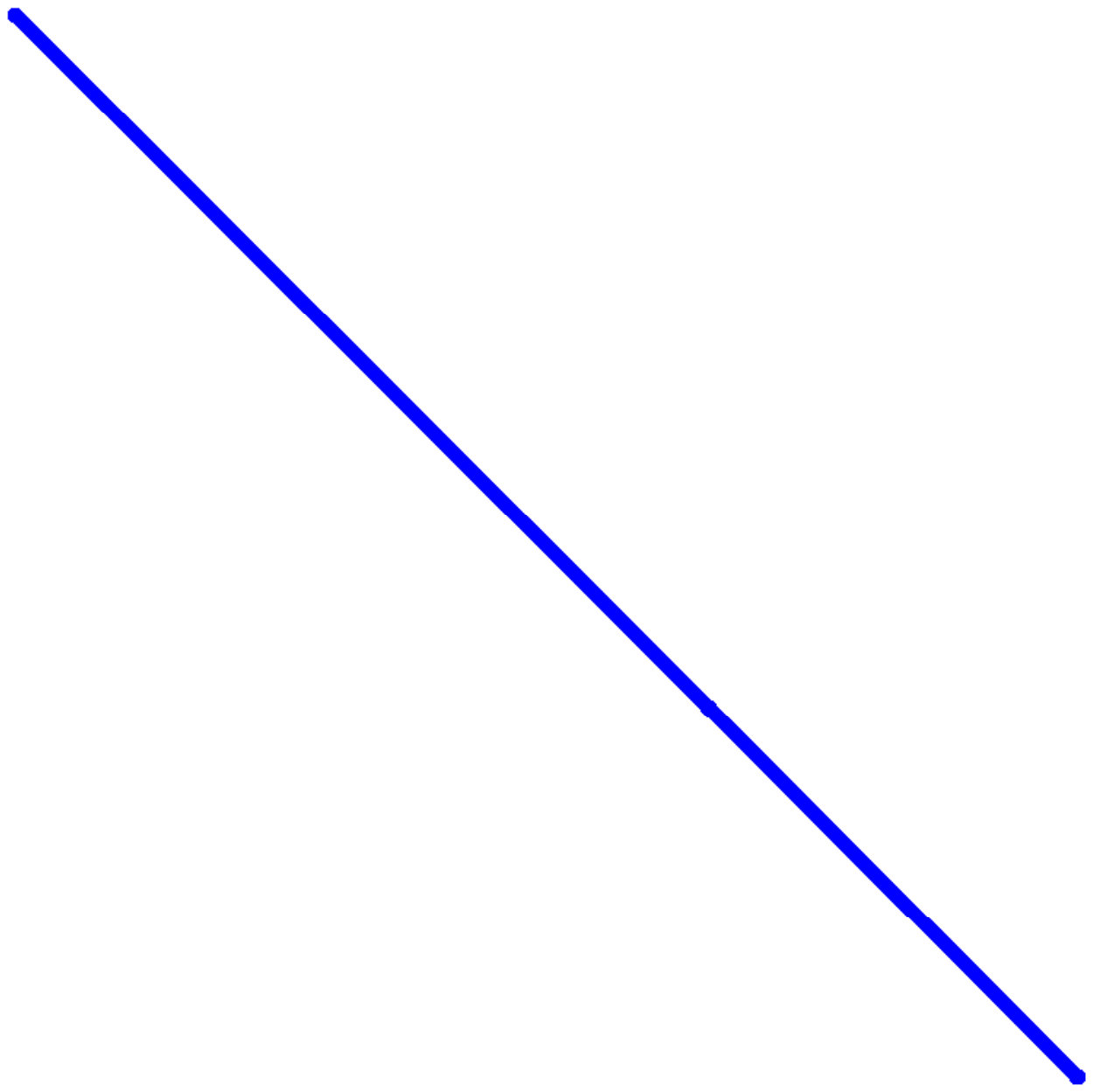}
\end{minipage}}
\hspace{0.3cm}
\subfloat[Smooth, level = 2]{
\begin{minipage}[c]{0.25\linewidth}
\centering
\includegraphics[trim=4cm 6cm 4cm 6cm,clip=true,width=\textwidth]{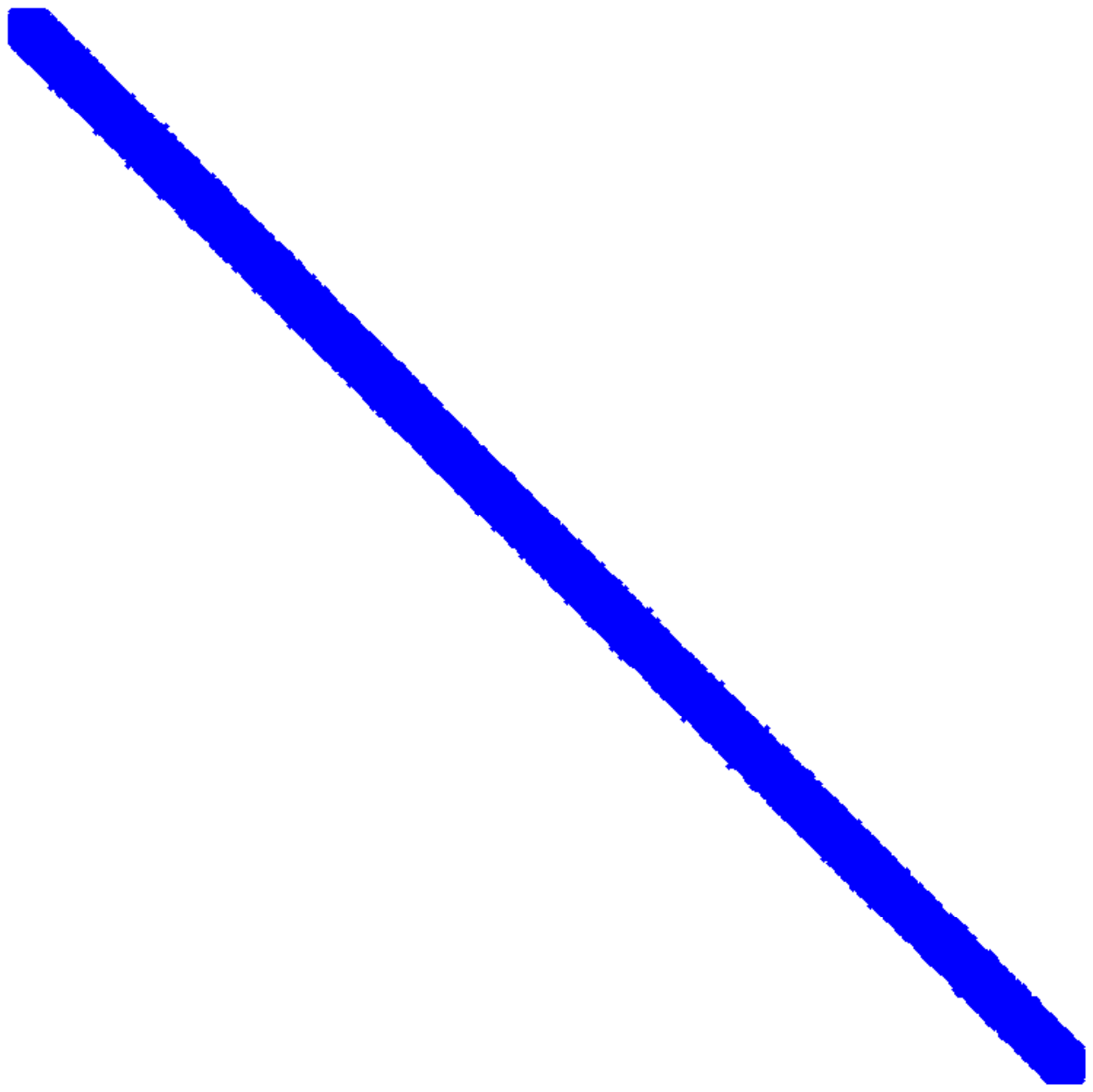}
\end{minipage}}
\hspace{0.3cm}
\subfloat[Smooth, level = 4]{
\begin{minipage}[c]{0.25\linewidth}
\centering
\includegraphics[trim=4cm 6cm 4cm 6cm,clip=true,width=\textwidth]{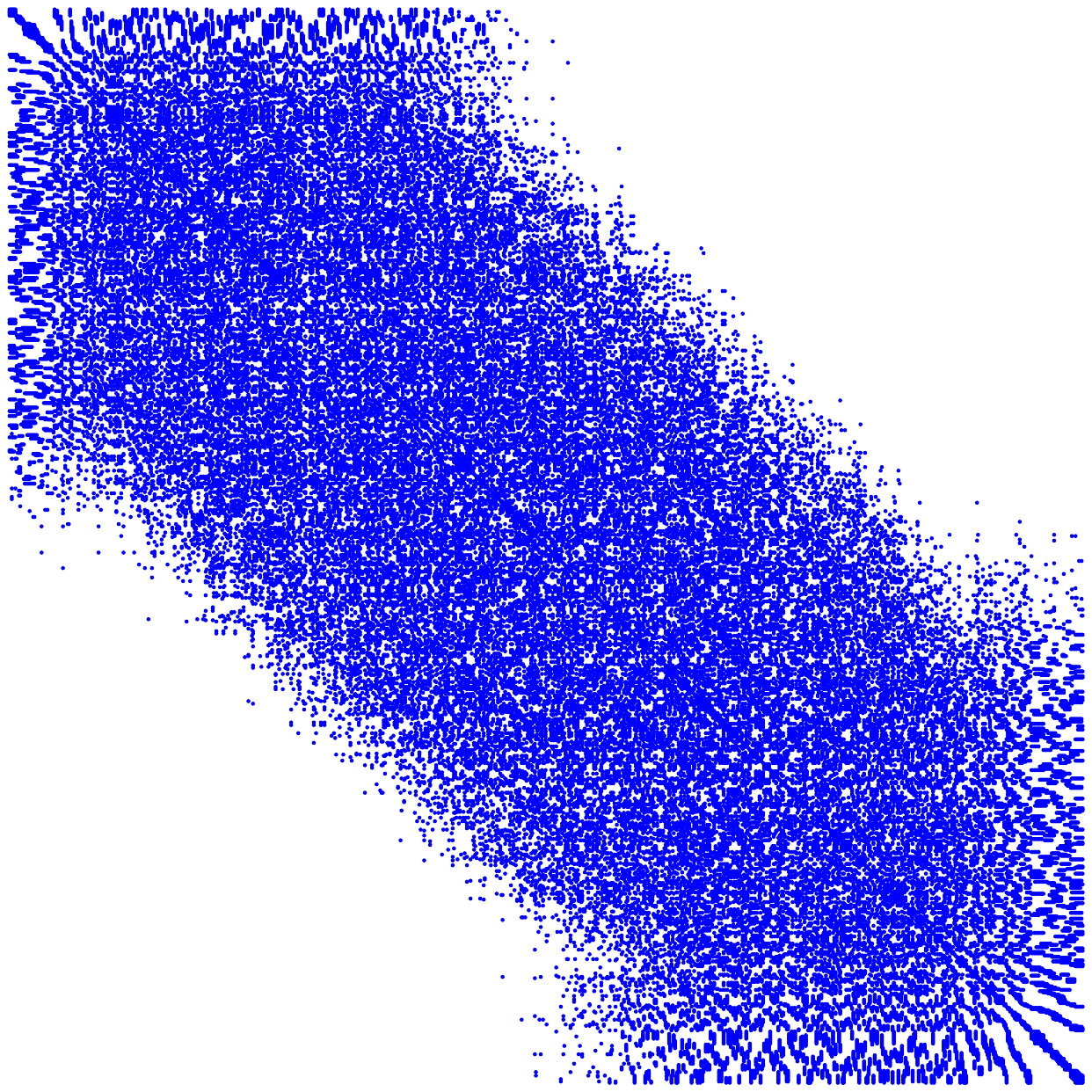}
\end{minipage}} \\
\subfloat[Unsmooth, level = 1]{
\begin{minipage}[c]{0.25\linewidth}
\centering
\includegraphics[trim=4cm 6cm 4cm 6cm,clip=true,width=\textwidth]{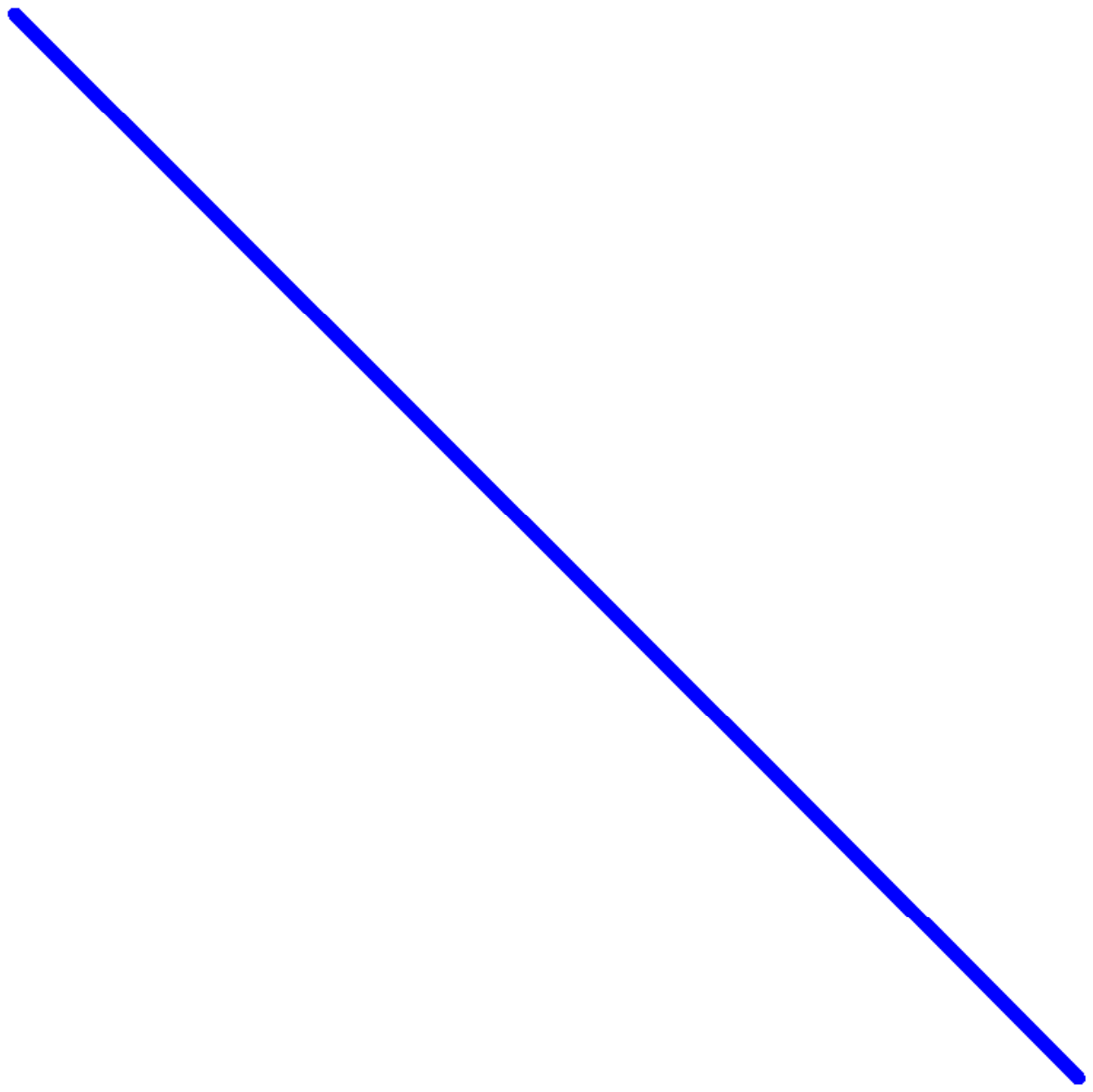}
\end{minipage}}
\hspace{0.3cm}
\subfloat[Unsmooth, level = 2]{
\begin{minipage}[c]{0.25\linewidth}
\centering
\includegraphics[trim=4cm 6cm 4cm 6cm,clip=true,width=\textwidth]{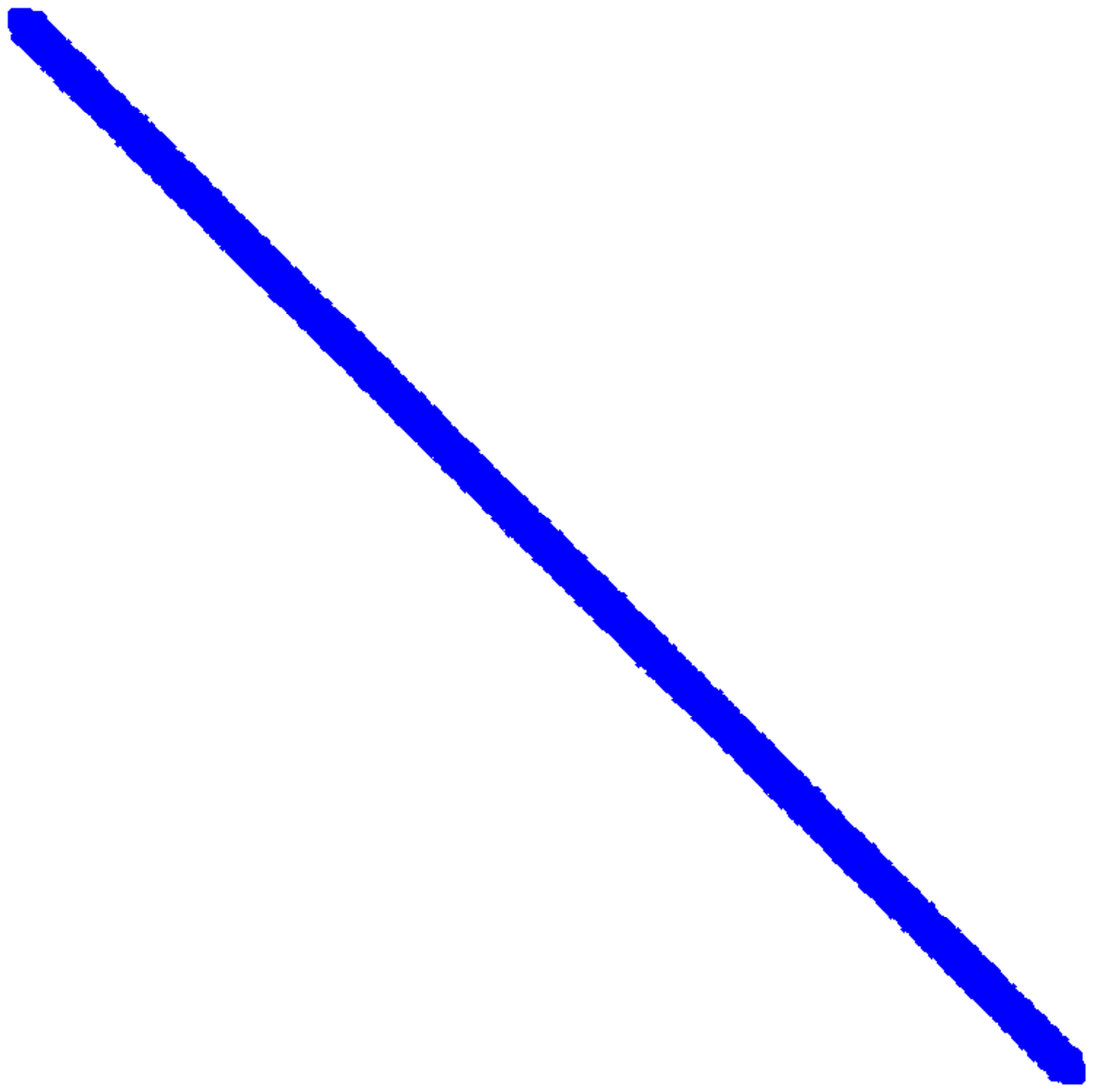}
\end{minipage}}
\hspace{0.3cm}
\subfloat[Unsmooth, level = 5]{
\begin{minipage}[c]{0.25\linewidth}
\centering
\includegraphics[trim=4cm 6cm 4cm 6cm,clip=true,width=\textwidth]{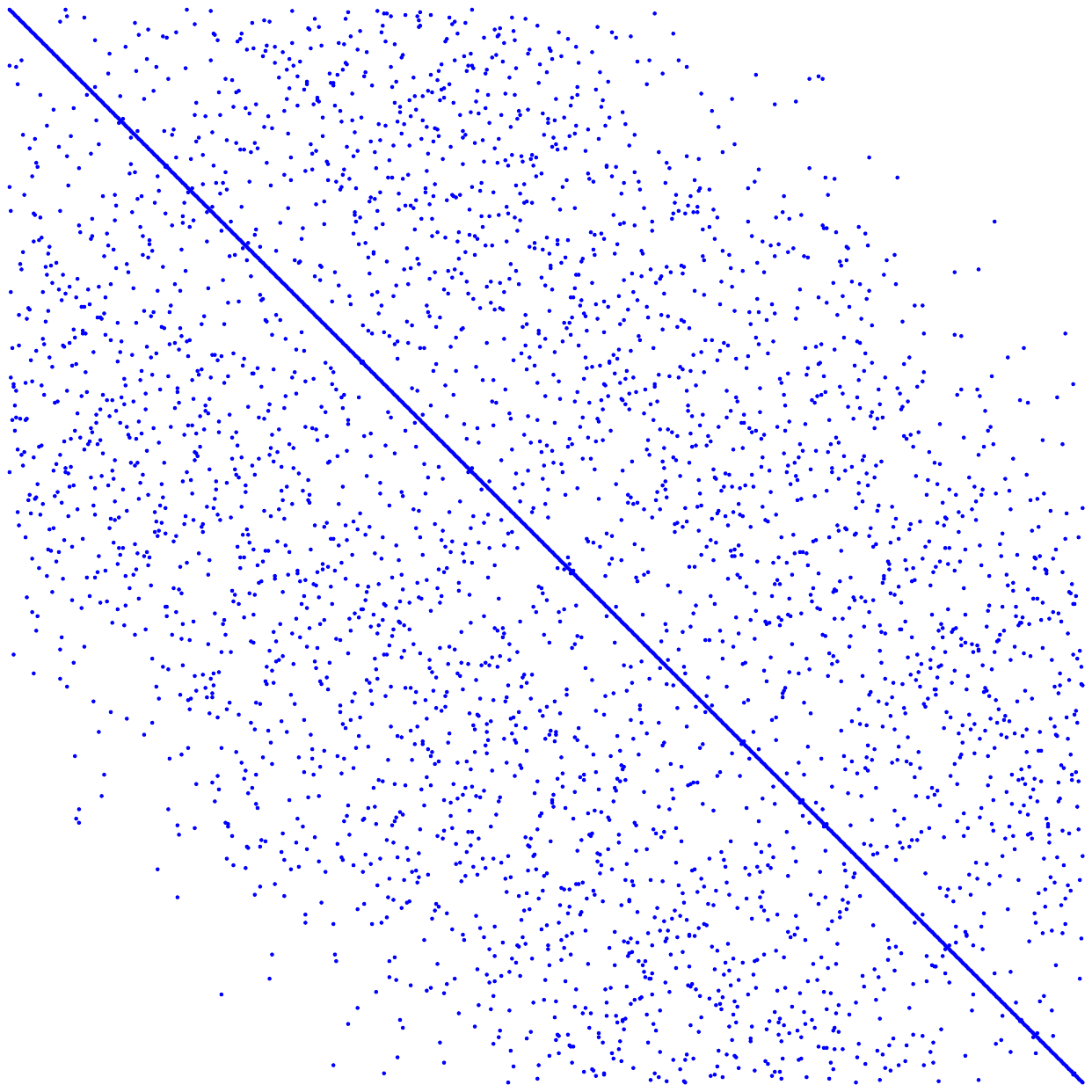}
\end{minipage}}
\caption{\emph{Comparison of sparsity pattern of CUSP smooth aggregation (top) with GAMPACK unsmooth aggregation (bottom) for 2D anisotropic Poisson problem on a $1000 \times 1000$ grid.}}
\label{fig:compare_sparsity}
\end{center}
\end{figure}

\noindent In Table (\ref{tab:sag_vs_agg}),  we compare the storage complexity of the hierarchy of smooth aggregation from \verb=CUSP= with unsmooth aggregation for a matrix related to a 2D anisotropic Poisson problem. For smooth aggregation, the density of the coarse grid matrices increase rapidly while it is almost fixed for unsmooth aggregation, ensuring low memory usage by unsmooth aggregation. In Figure (\ref{fig:compare_sparsity}) we compare the growth of density of coarse grid matrices for smooth aggregation  with the fixed sparsity pattern of unsmooth aggregation for the same problem. \\
\begin{table}
\begin{center}
\begin{tabular}{crrrcrrr}
\hline
& & & \verb=CUSP= & \verb=GAMPACK= & & & \\ \hline
level & unknowns &     nnz & nnz/row          & level & unknowns &     nnz & nnz/row \\ \hline
0     & 1000000  & 4996000 & \textbf{5.00}    & 0     & 1000000  & 4996000 & \textbf{5.00} \\
1     &  275264  & 4570716 & \textbf{16.60}   & 1     &  274979  & 1770339 & \textbf{6.44} \\
2     &   76002  & 4398710 & \textbf{57.88}   & 2     &   77433  &  525877 & \textbf{6.79} \\
3     &    8601  &  737445 & \textbf{85.74}   & 3     &   21823  &  147803 & \textbf{6.77} \\
4     &     923  &   82045 & \textbf{88.89}   & 4     &    6150  &   39298 & \textbf{6.39} \\
5     &     246  &   28140 & \textbf{114.39}  & 5     &     938  &    5128 & \textbf{5.47} \\
      &          &         &                  & 6     &     169  &     641 & \textbf{3.79} \\ \hline
\end{tabular}
\caption{\emph{Comparison of the hierarchy grid sparsity of CUSP smooth aggregation (left) and  GAMPACK unsmooth aggregation (right) for 2D anisotropic Poisson problem on $1000 \times 1000$ grid.}}
\label{tab:sag_vs_agg}
\end{center}
\end{table}

\begin{figure}
\begin{center}
\includegraphics[trim=0cm 4cm 0cm 4cm,clip=true,width=0.6\textwidth]{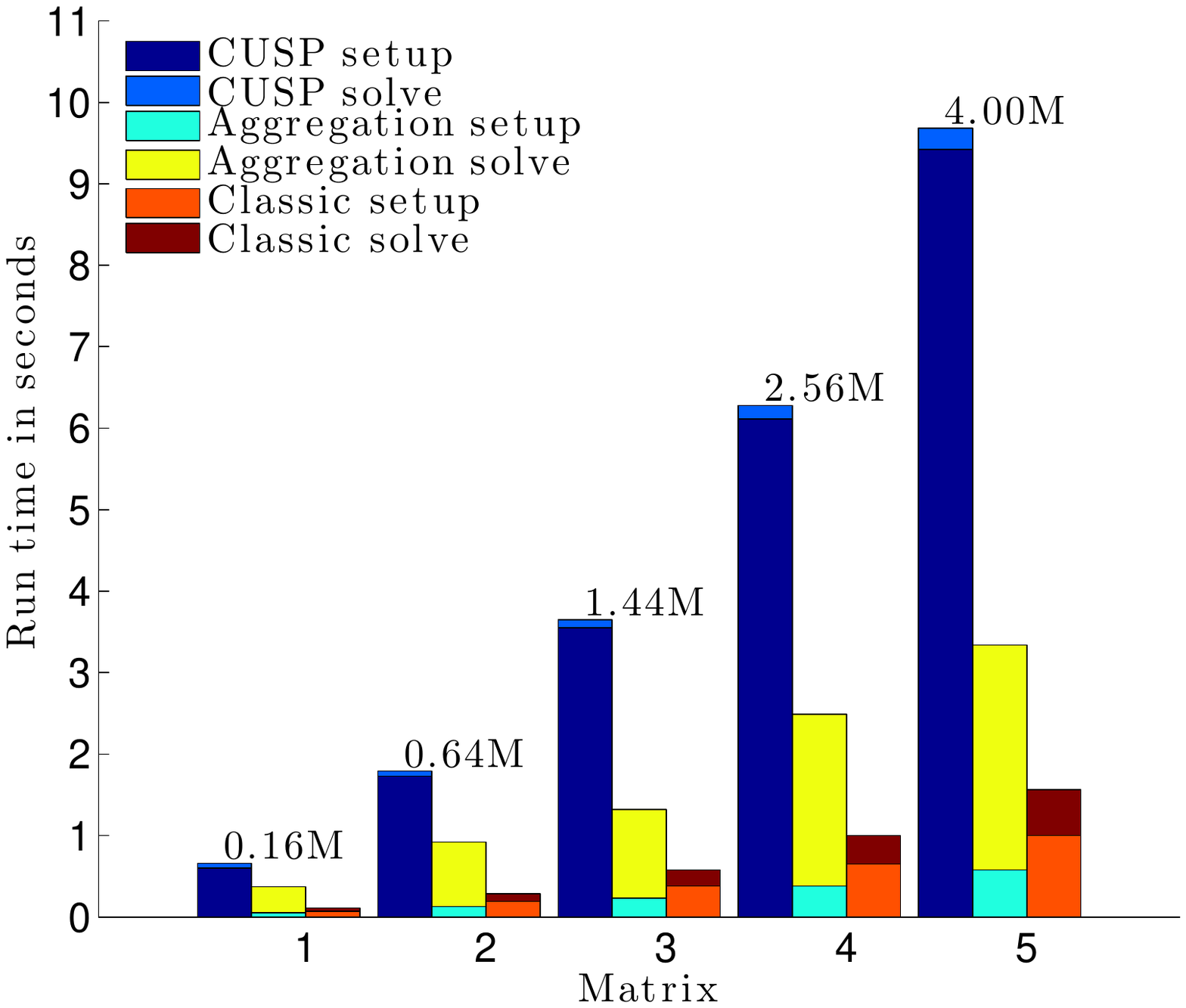}
\caption{\emph{Comparison of performance of CUSP smooth aggregation, GAMPACK aggregation, and GAMPACK classic AMG for 2D Poisson problem. Experiments ran on a single NVIDIA Tesla M2090 GPU. The numbers on the bars represent the size of the corresponding system in millions.}}
\label{fig:cusp_vs_agg}
\end{center}
\end{figure}

\noindent In Figure (\ref{fig:cusp_vs_agg}), we compare the setup and solve times with those of \verb=CUSP= smooth aggregation and \verb=GAMPACK= classical AMG for an anisotropic 2D Poisson problem. The setup times for \verb=CUSP= are significantly larger while the solve times are significantly smaller compared to that of \verb=GAMPACK=. Very small solve time of smooth aggregation is due to fewer number of iterations required for the convergence. This is because of very good interpolation between coarse  and fine grids. Since the setup is the dominant cost for smooth aggregation, the overall runtime for \verb=CUSP= smooth aggregation is also very large compared to \verb=GAMPACK= unsmooth aggregation. The large setup time for smooth aggregation is due to the density of the  coarse levels.   \\

\begin{figure}
\begin{center}
\includegraphics[trim=0cm 4cm 0cm 4cm,clip=true,width=0.6\textwidth]{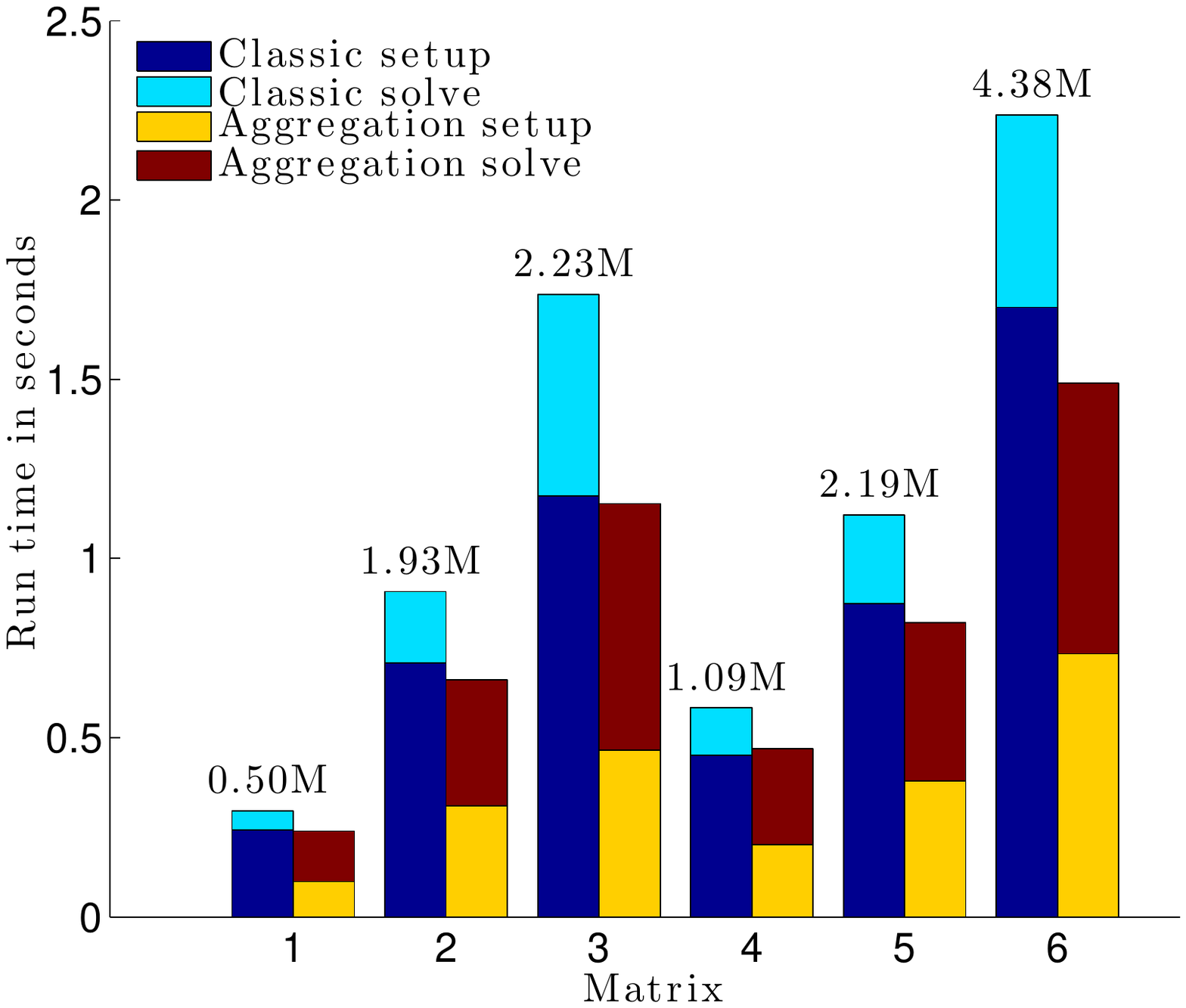}
\caption{\emph{Comparison of performance of GAMPACK classical AMG and aggregation AMG for RS and pressure matrices. Experiments ran on a single NVIDIA Tesla M2090 GPU. The numbers on the bars represent the size of the corresponding system in millions.}}
\label{fig:classic_vs_agg}
\end{center}
\end{figure}

\noindent In Figure (\ref{fig:classic_vs_agg}),  we compare the setup and solve times with classical AMG. Unlike classical AMG, the setup cost of aggregation is not dominant, and more than half of the time is spent on the solve. The results indicate that the solve stage of V-cycle classical AMG is about two times faster than aggregation AMG while the setup time of classical AMG is three times slower than that of aggregation AMG. Since a major component  of computation is  the setup for classical AMG, in overall, aggregation AMG is faster than classical AMG for the class of problems we considered.  For these test cases, \verb=CUSP= smooth aggregation could not construct the AMG hierarchy in a reasonable time, hence we do not compare the results with \verb=CUSP=.   \\

\noindent In order to study the scalability of the solver, we consider a series of down-scaled problems corresponding to the SPE10 pressure matrix problem \cite{christie2001tenth}. Since the larger matrices do not fit on one or two GPUs, we include the timings with multiple GPUs for both classical and aggregation AMG in Table (\ref{table:RS_matrices}). \\
\begin{table}
\begin{center}
\begin{tabular}{ccccc}
\hline
 N       & 1$\times$M2090 & 2$\times$M2090 & 3$\times$M2090 & 4$\times$M2090  \\ \hline
1.1M     & 0.469 (15)     & 0.387 (15)     & 0.430 (15)     & 0.409 (15) \\ 
         & 0.583 (10)     & 0.365 (10)     & 0.326 (10)     & 0.312 (11) \\ \hline 
2.2M     & 0.821 (15)     & 0.580 (15)     & 0.561 (15)     & 0.500 (15) \\
         & 1.121 (10)     & 0.684 (11)     & 0.526 (11)     & 0.454 (10) \\ \hline
4.4M     & 1.488 (15)     & 1.079 (15)     & 0.875 (15)     & 0.874 (15) \\
         & 2.237 (11)     & 1.302 (11)     & 0.953 (11)     & 0.774 (11) \\ \hline
8.9M     & -              & 1.636 (14)     & 1.345 (15)     & 1.177 (15) \\
         & -              & 2.556 (12)     & 1.856 (12)     & 1.564 (15) \\ \hline
13.5M    & -              & -              & 1.813 (15)     & 1.471 (14) \\
         & -              & -              & 2.769 (12)     & 2.159 (12) \\ \hline
\end{tabular}
\caption{\emph{Scaling : Timings in seconds for GAMPACK setup and solve for aggregation AMG, and classical AMG on single and multiple GPUs. The first column represents the number of cells in the original SPE10 pressure matrix problem. For each matrix, the first row gives the timings for aggregation AMG while the second row gives the timings for classical AMG. All timings are in seconds and include setup and solve to a relative tolerance of $10^{-6}$. The number of GMRES iterations are given in parentheses. }}
\label{table:RS_matrices}
\end{center}
\end{table}

\begin{figure}
\begin{center}
\includegraphics[trim=0cm 4cm 0cm 4cm,clip=true,width=0.6\textwidth]{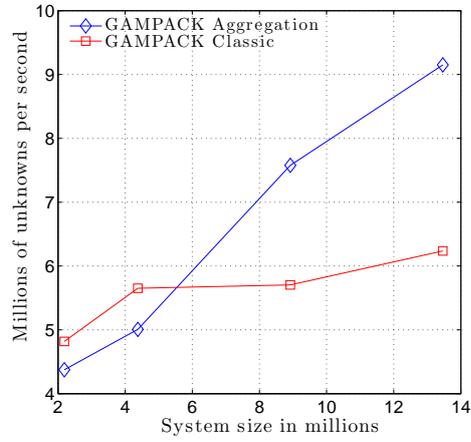}
\caption{\emph{Scaling of computation rate (in millions of unknowns per second) with system size. The rates include both setup and solve. The computations ran on 4 NVIDIA Tesla M2090 GPUs.}}
\label{fig:RS_scaling}
\end{center}
\end{figure}

\noindent In Figure (\ref{fig:RS_scaling}), we compare the performance of both solvers as the problem size increases. For the smallest problem ($\sim 2M$) we considered, classical AMG out-performs aggregation AMG, while aggregation AMG is about 50\% faster to classical AMG for the largest problem ($\sim 14M$) we considered. The performance of aggregation AMG\  improves at a faster rate compared to that of classical AMG.  \\

\begin{figure}
\begin{center}
\includegraphics[trim=0cm 4cm 0cm 4cm,clip=true,width=0.6\textwidth]{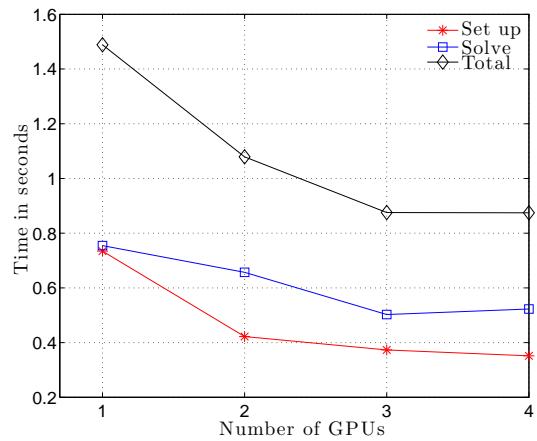}
\caption{\emph{Strong scaling : Computation time with number of GPUs. Computations ran for a system of four million unknowns on NVIDIA Tesla M2090 GPUs.}}
\label{fig:RS_strong_scaling}
\end{center}
\end{figure}

\noindent From Figure (\ref{fig:RS_strong_scaling}), we observe  significant decline in the performance of aggregation AMG as the number of GPUs used for the computation increases. This is because K-cycles require more number of coarse grid corrections  and GPUs are inefficient for systems of small size. For aggregation AMG to be efficient on multiple GPUs, the system size has to be sufficiently large. \\

%% file: conclusions.tex
\section{Conclusions and Future work}
\label{sec:conclusions}
We have presented a robust, efficient and scalable aggregation algebraic multigrid solver and have  verified the robustness with matrices from various applications. We compared its performance with a GPU accelerated classical AMG solver and a GPU\ accelerated smooth aggregation AMG solver, and observed that the setup cost of aggregation AMG is significantly low and the cost of the solution phase is high. In conclusion, aggregation AMG\ is extremely efficient for systems of sufficiently large size. Furthermore, for practical applications, that do not require high accuracy of the linear solver, aggregation AMG is extremely efficient\ due to very low setup costs. Our future work includes using MPI with multi-node computations  for very large systems and incorporating polynomial smoothers to improve the convergence of aggregation AMG.